\documentclass[12pt,a4paper]{article}
\usepackage[T1]{fontenc}
\usepackage[utf8]{inputenc}
\usepackage[english]{babel}
\usepackage{authblk}
\usepackage{amsmath,amssymb,amsthm}
\usepackage{graphicx}
\usepackage{color}
\usepackage{makecell}
\usepackage{hyperref}
\usepackage{amsfonts}
\usepackage{mathrsfs}
\usepackage{url}
\allowdisplaybreaks

\begin{document}
	
\begin{center}
	\textbf{{\Large Symmetric Branching Random Walks in Random Media: Comparing Theoretical and Numerical Results}}
\end{center}

\begin{center}
	\textit{{\large Kutsenko Vladimir,  Elena Yarovaya}}
\end{center}

\begin{abstract}

We consider a continuous-time branching random walk on a multidimensional lattice in a random branching medium. It is theoretically known that, in such branching random walks, large rare fluctuations of the medium may lead to anomalous properties of a particle field, e.g., such as intermittency. However, the time intervals on which this intermittency phenomenon can be observed are very difficult to estimate in practice. In this paper, branching media containing only a finite and non-finite number of branching sources are considered. The evolution of the mean number of particles with a random point perturbation and one initial ancestor particle at a lattice point is described by an appropriate Cauchy problem for the evolutionary operator. We review some the previous results about the long-time behavior of the medium-averaged moments $\langle m_{n}^{p}\rangle$, $p\geq 1$, $n\geq 1$ for the particle population at every lattice point as well as the total one over the lattice and present an algorithm for the simulation of branching random walks under various assumptions about the medium, including the medium randomness. The effects arising in random non-homogeneous and homogeneous media are then compared and illustrated by simulations based on the potential with Weibull-type upper tail. A wide range of models under different assumptions on a branching medium, a configuration of branching sources, and a lattice dimension were considered during the comparison. The simulation results indicate that intermittency can be observed in random media even over finite time intervals.

\end{abstract}


\section{Introduction}

We consider a continuous-time branching random walk (BRW) on the multidimensional lattice $\mathbb{Z}^\mathrm{d}$, $\mathrm{d}\in \mathbb{N}$, in a random and non-random  branching media. Branching random walks in non-random media have been extensively studied in a number of papers, see, e.g. \cite{ABY:1998, YarBRW:r, KY:19}, and the bibliography therein. Recall briefly some principal results related to BRWs in non-random media.

Let us start with an informal description of the process. Suppose that initially there is a single particle at a point $x\in{\mathbb Z}^\mathrm{d}$,
which then performs a random walk on ${\mathbb Z}^\mathrm{d}$ until reaching one of the
lattice point, {called the source of branching}, where it can  produce an arbitrary number of offspring or die.  Intensities of birth and death at the source are independent identically distributed  random variables (i.i.d.). We consider two cases of BRWs. The first of them, called {a spatially homogeneous} case, is in which every lattice point is a source of branching. The second one, called {a spatially non-homogeneous} case, is in which the random medium contains a finite number of sources. The underlying random walk is assumed to be symmetric, irreducible and homogeneous in space. All particles behave independently from each other and their prehistory.

A number of results were obtained for the model described above. The paper \cite{ABMY:2000} is currently contains the most complete and detailed analysis of BRWs in a spatially homogeneous random medium. The case of a spatially non-homogeneous random medium in which branching occurs at one lattice point, for instance, at the origin, is considered in~\cite{Yarovaya:2012}.  Moreover in~\cite{Yarovaya:2012} the results obtained earlier for the discrete Laplacian in the model of a BRW from~\cite{ABMY:2000} were extended to a wider class of an underlying symmetric random walk with a finite variance of random walk jumps.
The goal of these papers was to investigate the long-time behavior of the moments $\langle m_{n}^{p}\rangle$, $p\geq 1$, $n\geq 1$, for the particle population at every lattice point and the total one over the lattice, obtained by averaging over the medium
realizations.
The analysis of BRW in random media can be reduced to the analysis of an appropriate Cauchy problem for the evolutionary operator of a mean number of particles with a random point stationary perturbation, called {a random potential}, and one initial ancestor-particle at a lattice point as the initial condition. In particular,  the Feynman-Kac formula  for the solution of the non-homogeneous Cauchy problem is applied.
This approach make it possible to obtain  limit theorems for non-homogeneous BRWs in a random medium.

Large rare fluctuations of the random branching medium may lead to anomalous
properties of a particle field which is conventionally expressed in terms of the so-called intermittency concept~\cite{GM:90,GM:98,M:1994}. One of the first articles devoted to the intermittency problems in a  random media was developed in the paper by Zeldovich et al.~\cite{Zeld:87}. These problems were generalized and developed later in two fundamental papers by G{\"a}rtner and Molchanov \cite{GM:90,GM:98}. G{\"a}rtner and Molchanov rigorously formalized the concept of intermittency and developed some tools for studying the spatially homogeneous case. For the detailed discussion of the concept of intermittency of evolutionary random fields, and in particular the concept of ``strong centres'' of the field generation, see, e.g., \cite{ABMY:2000} and the bibliography therein.

The vast majority of the results obtained in the theory of BRWs are asymptotic in nature. At the same time, the study of BRWs at finite time intervals seems to be a difficult problem that has not yet been solved. On the other hand, in applied tasks, we can observe the behavior of the process only at a finite time interval. Thus, the main goal of this paper is to investigate whether it is possible to simulate qualitative and quantitative results predicted by the theory. In \cite{ermishkina} it is shown that such a simulation is possible for most theoretical results related to BRWs in a non-random medium. In particular, it is possible to show qualitatively different process behavior at different branching intensities and quantify the growth rates of the total and the local number of particles. Unfortunately, no similar studies, to the best of our knowledge, were conducted for a non-random medium. At the same time, the structure and behavior of the BRW in a random medium are more complex than in a non-random one. In particular, until now it is unclear whether it is possible to observe the phenomenon of intermittency at finite times, which is initially formulated in a strictly asymptotic formulation. Moreover, it is even unclear whether the phenomenon of intermittency allows modeling in random media at all. In addition, simulation of BRWs in random media requires dramatically more computer power than modeling BRW in a non-random media. In particular, in this work, we had to use a 48-core computer cluster for simulation.

The structure of the paper is as follows. In Section~\ref{S2}, we describe the model of BRW in random media in both homogeneous and non-homogeneous cases. In Section~\ref{S3}, we present the limit theorems and the definition of intermittency in the context of the BRW model. Finally, in Section~\ref{S4}, we present the algorithm and results of simulations of BRW with various initial conditions.

\section{Branching random walks in random media} \label{S2}

Assume that a branching process is formed at a lattice point $x\in \mathbb{Z}^\mathrm{d}$ by pairs
of non-negative random variables,
$\xi(x):=\xi(x,\omega):=(\xi^{-}(x,\omega),\xi^{+}(x,\omega))$, defined on a
probability space $(\Omega,\mathcal{F},\mathbb{P})$, where $\omega\in\Omega$ represent sample realizations of the field
$\xi(\cdot)$. The expectation with
respect to the probability measure $\mathbb{P}$ will be denoted by
angular brackets, $\langle\cdot\rangle$.
If the distribution
$\mathbb{P}$ of the field is invariant with respect to translations
$x\mapsto x+h$, for all $x,h\in\mathbb{Z}^\mathrm{d}$ (see, e.g.,
~\cite{ABMY:2000}), then the random
field $\xi$ is supposed {spatially homogeneous}, and $\Omega=(\mathbb{R}^2_+)^{\mathbb{Z}^\mathrm{d}}$.

In non-homogeneous case we consider a fundamentally different situation. A branching random medium is formed by a finite number of pairs of non-negative i.i.d. random variables $\xi(x_{i}):=\xi(x_{i},\omega)=(\xi^{-}(x_{i},\omega),\xi^{+}(x_{i},\omega))$, $i=1,2,\ldots, N$,
defined  on a probability space
$(\Omega,\mathcal{F},\mathbb{P})$ by the same law at the lattice points $x_{1},x_{2},\ldots, x_{N}$. It is assumed that
$\Omega=(\mathbb{R}^2_+)^{N}$, where $N<\infty$. In this case, the random environment is
spatially non-homogeneous, since the branching medium formed of
birth-and-death processes only at a few lattice points. In~\cite{Yarovaya:2012}  such model was investigated for $N=1$.

Let us describe the underlying random walk for both
cases. The random walk of particles is given by the infinitesimal
transition matrix $A=(a(x,y))_{x,y \in {\mathbb{Z}}^\mathrm{d}}$
and is assumed to be symmetric; $a(x,y)=a(y,x)$ for all $x,y \in
\mathbb{Z}^\mathrm{d}$; homogeneous in space: $a(x,y)=a(0,y-x)=a(y-x)$; irreducible
and regular: $\sum_{x \in {\mathbb{Z}}^\mathrm{d}}a(x)=0$ with $a(x)\geq
0$, $x\neq 0$, $a(0)<0$.  Since the random
walk is symmetric and homogeneous, the matrix $A$ satisfies
$\sum_{y\in \mathbb{Z}} a(x,y)=0$ and $\sum_{x\in \mathbb{Z}}
a(x,y)=0$, see, e.g., ~\cite{Yarovaya:2012}.
Under these conditions the transition probability $p(t,x,y)$ of the
random walk satisfies~\cite{Gih} the system of differential
Kolmogorov's backward equations
\begin{equation}\label{E:ptxy}
\partial_{t}p=\mathcal{A}p,\quad
p(0,x,y)=\delta_y(x),
\end{equation}
where the linear operator $\mathcal{A}$ (generated by the matrix $A$) defines the bounded linear operator $\mathcal{A}:\ell^p(\mathbb{Z}^\mathrm{d})\to\ell^p(\mathbb{Z}^\mathrm{d})$ for any $p\in[1, \infty]$ and acts with respect to the variable
$x$ as follows:
\[
\mathcal{A} p(t,x,y):=\sum_{x'}a(x,x')p(t,x',y),
\]
while
$y$ is treated as a parameter. Here, $\delta_y(x)$ is the
Kronecker symbol.

Additionally, we assume a finite variance of random walk jumps:
${\sum_{x \in {\mathbb{Z}}^\mathrm{d}} |x|^{2}a(x)<\infty}$. In particular, this class includes a simple symmetric random walk
defined by $a(x,y)= - a(0) / 2d$, for $|y -x| = 1, a(x, x) = a(0)$,
and $a(x,y) = 0$ otherwise. If we make the notation $\varkappa=- a(0)/2d$, then we can obtain
the simple symmetric random walk defined by an the operator $\varkappa
\Delta$ in the Cauchy problem~\eqref{E:ptxy}.

Therefore, if at time  $t=0$ there is one particle
 at a lattice point $x\in \mathbb{Z}^\mathrm{d}$, which is the source of branching, then  in the time
interval $[0, h)$, as $h \to 0$, it can either jump with the
probability $p(h,x,y)=a(x,y)h+o(h)=a(y-x)h+o(h)$ to the point
$y\neq x$, or produce  two particles (including itself) with the
probability $\xi^{+}(x)h+o(h)$, or die with the probability
$\xi^{-}(x)h+o(h)$, or survive at the point $x$ (no changes) with
the probability $1+a(0)h-(\xi^{+}(x)+\xi^{-}(x))h+o(h)$.

If a lattice point  $x\in \mathbb{Z}^\mathrm{d}$  is not  a source of branching, then in the time
interval $[0, h)$, as $h \to 0$, it the initial particle can only  jump with the
probability $p(h,x,y)=a(x,y)h+o(h)=a(y-x)h+o(h)$ to the point
$y\neq x$ or stay  at the point $x$ (no changes) with
the probability $1+a(0)h+o(h)$.

In both models, evolution of the system of
particles on ${\mathbb{Z}}^\mathrm{d}$ is described by the number
$\mu_{t,\omega}(y)$ of particles at time $t$ at each point $y\in
{\mathbb{Z}}^\mathrm{d}$  and the total particle population size
$\mu_{t,\omega}:=\sum_{y\in \mathbb{Z}^\mathrm{d}}\mu_{t,\omega}(y)$ under the
assumption that at time $t=0$ the system consists of one particle
at a point $x$, that is, $\mu_{0}(y)=\delta_{y}(x)$. It is also
assumed that evolution of newborn particles obeys the same law
independent of other particles and  the prehistory.

For every $n \in \mathbb{N}$ the ``quenched''  moments, see, e.g.,~\cite{ABMY:2000}, are introduced as:
\begin{align*}m_n^p(t,x,y):&=m_n^p(t,x,y, \omega)  = \left( \mathbb{E}_x \mu_{t,\omega}^n\left(y\right)\right)^p,\\
m_n^p(t,x):&= m_n^p(t,x,\omega) = \left(\mathbb{E}_x \mu_{t,\omega}^n \right)^p.
\end{align*}
Here, the label $\omega$ is referred to the fixed realization of the
branching medium, and the subscript $x$ indicates the initial
position of the initial particle.
We denote by $\langle \cdot \rangle$ the expectation with respect to the medium probability measure. For every $n \in \mathbb{N}$ the ``annealed'' moments are determined as $\langle m^p_n(t,x,y)\rangle $ and $\langle m^p_n(t,x)\rangle$, respectively.

Finally, we briefly describe a BRW in a non-random media. For a detailed description, see, e.g., \cite{ABY:1998, YarBRW:r, BMY:20}. In this case, the branching intensities are assumed to be constant.The rest of the description remains the same with only one exception. In a non-random media, there is no concept of a quenched moment. The moments of the number of particles are non-random variables and are defined as

\begin{align*}m_n^p(t,x,y)  = \left( \mathbb{E}_x \mu_{t}^n\left(y\right)\right)^p,\\
 m_n^p(t,x) = \left(\mathbb{E}_x \mu_{t}^n \right)^p.
\end{align*}

\section{Asymptotic results}\label{S3}

In this section, we recall the main results, based on the papers \cite{ABMY:2000} and  \cite{Yarovaya:2012}. The generating functions associated with the random
variables $\mu_{t,\omega}(y)$ and $\mu_{t,\omega}$ are defined for $z>0$ by
\[
F(z;t,x,y):=\mathbb{E}_{x}(e^{-z\mu_{t,\omega}(y)}),\quad
  F(z;t,x):=\mathbb{E}_{x}(e^{-z\mu_{t,\omega}}).
\]
Functions
$F(z;t,x,y)$, $F(z;t,x)$ satisfy  the equation for BRWs in {a spatially homogeneous medium}
\begin{equation}\label{E:Fh}
\partial_{t}F=\mathcal{A} F-
(\xi^{+}(x)+\xi^{-}(x))F+\xi^{+}(x)F^{2},
+\xi^{-}(x)
\end{equation}
and in {a spatially non-homogeneous medium}
\begin{equation}\label{E:Fnonh}
\partial_{t}F=\mathcal{A} F-
\sum_{i=1}^{N}\Delta_{x_{i}}\left((\xi^{+}(x_{i})+\xi^{-}(x_{i}))F+\xi^{+}(x_{i})F^{2}
+\xi^{-}(x_{i})\right),
\end{equation}
where $\varDelta_{x_{i}}=\delta_{x_{i}}\delta_{x_{i}}^{T}$, $i=1,2,\ldots,N$,
with the same initial conditions for~\eqref{E:Fh} and~\eqref{E:Fnonh} $F(z;0,x,y)=e^{-z\delta_y(x)}$ and
$F(z;0,x)=e^{-z}$, respectively, see, e.g., ~\cite{ABMY:2000,Yarovaya:2012}.

Let us introduce the concept of random potential as in \cite{ABMY:2000}:
\[
V(x):=\xi^{+}(x)-\xi^{-}(x), \quad
x\in\mathbb{Z}^\mathrm{d}.
\]

Then the moment functions for homogeneous medium $m_{n}(t, x, y), m_{n}(t, x)$ satisfy the chain of linear differential equations
$$\partial_{t} m_{n}=\mathcal{A} m_{n}+V(x) m_{n}+\xi^{+}(x) g_{n}\left(m_{1}, \ldots, m_{n-1}\right)$$
with the initial conditions $m_{n}(0, \cdot, y)=\delta_{y}(\cdot)$ and $m_{n}(0, \cdot) \equiv 1$, where $g_{1} \equiv 0$ and for $n \geq 2$:
\[
g_{n}\left(m_{1}, \ldots, m_{n-1}\right):=\sum_{i=1}^{n-1}\left(\begin{array}{c}
n \\
i
\end{array}\right) m_{i} m_{n-i}
\]

The same equations for a homogeneous medium have the following form:
$$\partial_{t} m_{n}=\mathcal{A} m_{n}+V(0) m_{n}+\xi^{+}(0)\delta_0 (x) g_{n}\left(m_{1}, \ldots, m_{n-1}\right)$$
with the initial conditions as in non-homogeneous case.

 In particular, the first-order moments  for homogeneous medium satisfy the following equation
$$
\partial_{t} m_{1}=\mathcal{A} m_{1}+V(x) m_{1}$$

In turn, the first-order moments  for non-homogeneous medium satisfy the following equation
\[\partial_{t} m_{1}=\mathcal{A} m_{1}+V(0)\Delta_0 m_{1}.\]

Apparently, for the first time, the limit theorems for moments of a simple symmetric BRW with a source of branching at each lattice point in a random medium were appeared in~\cite{ABMY:2000}. In~\cite{ABMY:2000} the intensities of splitting and death were assumed to be i.i.d for all sources. The proof of limit theorems from equations for moments is mainly based on the Feynman-Kac representations. To apply this technique, the following restriction on the random potential was introduced:
\begin{equation} \label{cond1}
    {\left\langle\left(\frac{V^{+}(0)}{\ln _{+} V^{+}(0)}\right)^\mathrm{d}\right\rangle<\infty},
\end{equation}
where $\ln_{+}(x):=\max\left(\ln(x),1\right)$.

Under these assumptions, limit theorems for annealed moments were proved. In paper \cite{Yarovaya:2012}, these theorems were generalized for an arbitrary symmetric, irreducible, and homogeneous in space random walk with a finite variance of walk jumps and a non-homogeneous random medium. For this purpose, a similar condition for the branching potential was introduced as follows:
\begin{equation} \label{cond2}
\lim\limits_{t \to \infty} \frac{t}{\ln \left\langle {\rm e}^{p V(0)t}\right\rangle} = 0, \quad p> 0.
\end{equation}

Under assumptions~\eqref{cond1}  and \eqref{cond2} the following results for the  annealed moments $\langle m_{1}^{p}\rangle$ and $\langle m_{n}^{p}\rangle$  are valid regardless of the homogeneity of the medium:

\begin{equation}
\label{theorem1}
\lim _{t \to \infty} \frac{\ln \left\langle m_{1}^{p}\right\rangle}{\ln \left\langle {\rm e}^{pVt}\right\rangle}=1,
\end{equation}

\begin{equation}
\label{theorem2}
\lim _{t \to \infty} \frac{\ln \left\langle m_{n}^{p}\right\rangle}{\ln \left\langle {\rm e}^{p nVt}\right\rangle}=1.
\end{equation}

The condition~\eqref{cond2} is valid for random
potentials having, e.g., the Weibull-type and the Gumbel-type upper tails~\cite{Yarovaya:2012}. Note that condition~\eqref{cond2} cannot be fulfilled if the potential is non-random. Namely, for the first moment of the supercritical BRW in a non-homogeneous and non-random medium condition~\eqref{cond2} takes the following form:
\[\lim _{t \to \infty} \frac{t}{\ln{m_1}} = \frac{t}{\ln{{\rm e^{\lambda t}}} } = \frac 1\lambda > 0, \quad \lambda>0. \]

For additional information about the Feynman-Kac representation and the proof of limit theorems, see \cite{ABMY:2000, Yarovaya:2012,GM:90}.

For studying the field of ``quenched''  moments we will use commonly accepted notation. Let $\left\{\eta(t, x) ; x \in \mathbb{Z}^\mathrm{d}\right\}$ be a family of non-negative spatially homogeneous random fields and \[\Lambda_{p}(t)=\ln \left\langle\eta(t, 0)^{p}\right\rangle < \infty, \quad t \geq 0,~p \in \mathbb{N}.\]
Given two functions $f, g: \mathbb{R}_{+} \to \mathbb{R}$, we will write $f \ll g$ if $g(t)-f(t) \to \infty$ as $t \to \infty$.
The random fields $\left\{\eta(t, x) ; x \in \mathbb{Z}^\mathrm{d}\right\}$ are called intermittent \cite{GM:90} if they are ergodic and
\begin{equation}
\label{LAM}
\Lambda_{1}(t) \ll \frac{\Lambda_{2}(t)}{2} \ll \frac{\Lambda_{3}(t)}{3} \ll \ldots
\end{equation}
Note that if so-called Lyapunov exponents
$
\lambda_{p}=\lim\limits_{t \to \infty} \frac{\Lambda_{p}(t)}{t^\beta}$, $p \in \mathbb{N}$,
exist for some $\beta \geqslant 1$, and
$
\lambda_{1}<\frac{\lambda_{2}}{2}<\frac{\lambda_{3}}{3}<\ldots \text { , }
$
then condition~\eqref{cond1} is valid.

The intermittency means that as $t \to \infty$ the main contribution to each moment function is carried by higher and higher and more and more widely spaced ``peaks'' of the random field.  Let us consider the field $m_n^p$ we are studying. The concept of intermittency applied to $m_n^p$ implies that
$\langle m_{n}(t) \rangle^p \ll \langle m_n^p(t) \rangle,\text{ as } t \to \infty $, particularly $\langle m_{1}(t) \rangle^2 \ll \langle m_1^2(t) \rangle$. This condition shows that the growth rate of the moment $\langle m_n^p(t) \rangle$ is anomalous as compared, e.g., with Gaussian. Namely, $\langle m_n^p(t) \rangle$ growth ``irregularly'' with its number $p$: for large $t$ the second moment is much larger than the square of the first moment, the fourth is much larger than the product of the first and the third and the square of the second and so on \cite{GM:90}.

Large rare fluctuations of a random environment may lead to anomalous properties of a particle field such as intermittency. Namely, within intermittency the main contribution to each moment function is carried by high and rare peaks of underlying random field. In our case the main contribution to $\langle m_n^p(t) \rangle$ is carried by such moments $m_n^p(t,\omega)$ which are based on ``lucky'' medium realizations $\omega$. For example, with supercritical intensities on an area around the position of the initial particle at $t=0$. This interpretation  is not obvious and is given in~\cite{GM:90}.

Let for simplicity consider $\left\langle m_{1}^{p}\right\rangle$, and assume that potential has the Weibull-type upper tail $ \mathbf{P}(V> r) \sim {\rm e}^{-cr^{-\alpha}}, \quad r\to \infty$, where $\alpha>1$, $c>0$.
In this case the representation~\eqref{theorem1} takes the following form
 \[\lim _{t \to \infty} \frac{\ln \left\langle m_{1}^{p}\right\rangle}{t^{\frac{\alpha}{\alpha-1}}}=C\cdot p^{\frac{\alpha}{\alpha-1}}.\]

Thus, for the Lyapunov exponents of the field $m_1^p$ we have \[\frac{\lambda_{p}}{p} = \frac{1}{p}\lim\limits_{t \to \infty}  \frac{\ln \left\langle m_{1}^{p}\right\rangle}{t^{\frac{\alpha}{\alpha-1}}}= C\cdot p^{\frac1{\alpha-1}}.\]

Therefore, the sequence of $\frac{\lambda_p}{p}$ is strictly increasing as function of $p$, which, in turn, implies the intermittency of $\left\langle m_{1}^{p}\right\rangle$. The same result is valid for $\left\langle m_{n}^{p}\right\rangle$, where $n \in \mathbb{N}$.

\section{Numerical Results} \label{S4}

The vast majority of the results obtained in the theory of branching random walks are asymptotic in nature~\cite{ABMY:2000}.
At the same time, the study of BRW at finite time intervals seems to be a difficult task that has not yet been solved satisfactory enough and remains actual. The main problem here is that, in practice, we are able to observe only the behavior of the process at a finite time interval. Thus, the main goal of the simulation in the paper is to investigate whether it is possible to obtain qualitative and quantitative results predicted by the theory already at finite times and how to do it in practice. A similar task was considered in \cite{ermishkina} for BRWs in  non-random media where the growth rates of the total and the local number of particles were calculated.
However, to the best of our knowledge there are no similar studies related to the simulation of BRWs in random media. This may be explained by the fact that  the structure and behavior of the BRW in a random media are more complex than in a non-random one. Here, the phenomenon of intermittency is of the most significant interest and importance. This phenomenon can manifest itself so sharply that it will not be possible to simulate the BRW for a sufficiently long time using reasonable computing power. Moreover, until now it was unclear whether it is possible to observe the intermittency phenomenon at finite times. It is worth noting also that branching processes are often used in practice, see e.g. \cite{GMS-B:10, MKY:21, MPR:2007}. At the same time, considering the random spatial structure allows a strong generalization of the existing models.

\subsection{Description of the simulation algorithm}

For the simulation, we used the R v3.6.3 data analysis environment by applying the parallel programming techniques on a computer cluster with $48$ CPU cores. The algorithm for modeling the BRW was based on the decomposition of a BRW into a set of exponential and polynomial random variables. Similar algorithms were used in \cite{ermishkina} to simulate BRWs in a non-random media.

First of all, let us specify the characteristics for the BRW under consideration involved in the simulation. The length of the time interval $[0, T_{max}]$ on which the simulation is performed is fixed first. Next, we specify the dimension $\mathrm{d}$ of the integer lattice. In the simulation we can consider a lattice containing a potentially infinite number of points. To do this, one can use dynamic memory allocation for an array in which information about the lattice is recorded. It is also possible to work with a finite sufficiently large area of the lattice, beyond which the walk does not go out for the time $T_{max}$. In this case, it is enough to use a bounded array containing the part of the lattice under consideration. In our simulation, we decided to apply the second approach and used a cubic lattice area with the side of $100$ vertices. This finite lattice is denoted below by $\mathcal{Z}^\mathrm{d}$. After defining the dimension and the size of the integer lattice, it is necessary to set the intensities of splitting and death of particles at each point of the lattice. These values are denoted as $\tilde \xi_+(x)$ and $\tilde \xi_-(x)$, respectively. Outside of the sources, the quantities $\tilde \xi_+(x)$ and $\tilde \xi_-(x)$ are assumed to be zero. At the sources, these values are set depending on the model type, e.g., as sample realizations of random variables. Finally, one needs to select the type of random walk. In this paper, we decided to use a simple symmetric random walk. Within such framework, the particle waits for an exponentially distributed time $Exp(\varkappa)$ from the moment of birth and then jumps equally likely to one of the neighboring vertices of the lattice. BRW with such a walk is easier to model than a BRW with more general random walk. At the same time, BRW with a simple symmetric random walk has a full list of effects of interest what is demonstrated in~\cite{ABMY:2000}.
So, for simplicity, we further describe the BRW simulation algorithm for a simple symmetric random walk.

Let us define the state of the BRW as a set of tuples $(k,t_b (k), t_d (k), x (k))$. Each tuple describes a particle marked by a unique number $k$, which  is born at the point $x \in \mathcal{Z}^\mathrm{d}$ at the moment $t_b (k)$ and ``evolves'' at the moment $t_d (k)$. During the time interval $[t_b(k), t_d(k)]$ this particle is located at the point $x (k)$. By ``evolution'' we mean a jump to another point, splitting, or the death of a particle. For convenience, we consider the transition of a particle to another vertex as the death of a particle at the starting point and the birth of a new particle at the destination point.

In all the simulations under consideration, at the initial moment of time, there is a single point located at zero. In this case the initial state of the BRW consists of a single element $(1, t_b(1), t_d(1), x (1)) = \left(1,0,\text{\it NA}, (0,\dots,0)\right)$. With the symbol {\it``NA'' (not available)} we denote the times of death that are not yet known. If there is at least one particle with $t_d = \text{\it NA}$, then the algorithm performs a step. At each step of the algorithm, among all particles with $t_d = \text{\it NA}$ the particle with the smallest $t_b$ is selected and processed. Let this particle be described by the vector $\left (k, t_b (k),\text{\it NA}, x(k)\right)$. The selected particle is processed according to the following rules:

\begin{enumerate}
    \item To $t_d (k)$ the value $t_b (k) + Exp (\varkappa)$ is assigned.
    \item  One of three types of evolution: ``jump'', ``split'' and ``die'' is selected with the probability $(\frac{\varkappa}D, \frac{\tilde \xi_+}D, \frac{\tilde \xi_-}D)$, respectively, where $D = \varkappa+ \tilde \xi_+\left(x(k)\right) + \tilde \xi_-\left(x(k)\right) $.
    \item If ``jump'' is chosen at step 2), then equally likely one of the neighbouring vertices next to $x (k)$ denoted by $x (k+1)$ is selected. A new particle $\left (n+1, t_d(k),\text{\it NA}, x(k+1)\right)$ is added to the set of tuples.
    \item If ``split'' is chosen at step 2), then two particles  $\left (n+1, t_d(k),\text{\it NA}, x(k)\right)$ and $\left(n+2, t_d(k),\text{\it NA}, x(k)\right)$ are added to the set of tuples. $\left(n+1,t_d(k),\text{\it NA}, x(k)\right)$ and $\left(n+2,t_d(k),\text{\it NA}, x(k)\right)$.
    \item If ``die'' is chosen in step 2), then nothing is added to the set of tuples.
\end{enumerate}

The algorithm stops working if there are no particles left with $t_d = \text{\it NA}$ or if all particles have $t_d > T_{max}$. In addition, the algorithm is interrupted if there are totally more than $1000$ particles on the lattice. We use the stopping rule since the data processing time increases exponentially with the exponential increase of the amount of particles. The time of such an interruption of the algorithm is denoted by $T_{stop}$. $T_{stop}$ can be significantly less than $T_{max}$, especially in the case of random media, when the intensity of splitting can be much higher than the intensity of death.

In this paper, we study only the total number of particles on the lattice at time $T$, denoted as $\mu(t)$. It is defined as the number of particles existing at time $t$, i.e. $\mu(t) = \sum_k\mathbb{I} \left\{t \in [t_b (k), t_d (k)]\right\}$, where $\mathbb{I}$ is the indicator function. In the case of fixed medium intensities, if this value of $\mu(t)$ experiences a sharp increase, then it is exponential \cite{ABY:1998}. Thus, we assumed that if the number of particles reaches $1000$, then $\mu (t)$ experiences exponential growth, and the influence of randomness can be neglected. To predict the number of particles after $T_{stop}$, we used a generalized linear regression with the exponent as a link function. The regression model was formulated as follows

$$ \quad \log \mu(t) = \alpha t +\varepsilon_t, $$
where $\varepsilon_t$ are i.d.d. normal variables with an equal variance. For model fitting, we used only $t>T_{100}$ where $T_{100}$ is the time of the first appearance of a hundred particles on the lattice. The total number of particles on the segment $T_{stop},T_{max}$ was predicted based on the fitted model.

The assumption that the exponential regression is a good approximation for the total number of particles needs to be verified.
To test this we used $250\,000$ simulated BRW trajectories in a homogeneous non-random medium with intensities $ \tilde \xi_ + = 2 $ and $ \tilde \xi_ - = 1$. From these trajectories, only those where the total number of particles reached $1000$ were selected. The actual values of $ \mu(t) $ were compared with the regression prediction on the time interval $ [T_ {100}, T_{stop}]$ and the coefficient of determination was calculated.
The average coefficient of determination was equal to $0.996$; the lowest was equal to $0.939$. The average error of the number of particles was equal to $0.08$; the maximum absolute error was equal to $ 13.0 $ of particles. We decided that this result indicated a decent approximation quality. Therefore we used exponential regression as part of the $\mu(t)$ simulation algorithm.

The algorithm described above generates the value $\mu (t)$ on the time interval $[0, T_{max}]$ for some pre-fixed initial conditions. Consider $M$ runs of the simulation algorithm with the results $\mu (t)_1, \dots, \mu(t)_M$. In order to evaluate $m_n (t) = \mathbb{E}_p\mu (t)$ , we can apply the Monte Carlo method, namely, put:

$$\hat m_n(t) = \frac 1{M} \sum\limits_{i = 1}^{M}  \mu_i^n(t) ,$$
This value converges to $m_n(t)$ according to the law of large numbers as $M \to \infty$. In our work $M =  1000$.

In the case of a random media, the value of $m_n (t)$ is random, which we emphasize by denoting it as $m_n(t,\omega)$. The sample realization of a random media is the sample realization of the intensities $\tilde \xi_+(x, \omega), \tilde \xi_-(x, \omega)$. Therefore, for a fixed $\omega$ quenched moment $m_n(t,\omega)$ can be estimated by the Monte Carlo method as

$$\hat m_n(t,\omega) = \frac 1{M} \sum\limits_{i = 1}^{M}  \mu_i^n(t,\omega) ,$$

Let the values $\hat m_n (t, \omega_1), \dots, \hat m_n(t,\omega_1), \dots, \hat m_n(t,\omega_{M_1})$ be estimated for various $ (\omega_1, \dots, \omega_{M_1})$, respectively. Then the estimate of the annealed moments is $\langle m_n^p(t) \rangle$ can be calculated as

$$\widehat { \langle  m^{p}_n(t) \rangle} = \frac 1{M_1} \sum\limits_{k = 1}^{M_1} \hat m_n^p(t,\omega_k)$$

In our work, $M_1$ is chosen to be equal to $250$. We assumed that if the effects are not sufficiently noticeable, we can increase the sample of $\widehat { \langle  m^{p}_n(t) \rangle}$ by increasing $M_1$, without reevaluation of previously estimated $\widehat { \langle  m^{p}_n(t) \rangle}$. However, it turned out that a sample size of $250$ is enough to observe the effects of interest.

We need to emphasize that in a non-random medium, the value of $m_n(t)$ is not random. Therefore, formally, the Monte Carlo method can be stopped at first step and $M$ simulations. However, for the convenience of comparing random and non-random media, we used the second averaging as in the random case. For the BRW in a non-random medium, we define the ``pseudo-annealed'' moment $\widehat { \left[ m^{p}_n(t) \right]}$ as follows:

\[\widehat { \left[   m^{p}_n(t) \right] } = \frac 1{M_1} \sum\limits_{k = 1}^{M_1} \hat m_n^p(t).\]

Note that $\widehat{\left[m^{p}_n(t) \right]}$ and $\hat m^{p}_n(t)$ are normally distributed, $\widehat{\left[m^{p}_n(t) \right]} \xrightarrow{\text{a.s}} m_n^p$ and ${\hat m^{p}_n(t) \xrightarrow{\text{a.s}} m_n^p}$ as $M \to \infty$ and $M_1 \to \infty$. In particular, if the hypothesis about the normality of the value $ \hat m^{p}_n (t)$ is not rejected, then the Monte Carlo method is implemented correctly.

\subsection{Models under consideration}

Recall that we assume that the random walk is simple and symmetric in all models. In addition, the intensity of the walk is assumed to be equal to $1$. That is, the particle waits for an $Exp(1)$ time and then moves equally likely to one of the neighboring points.

We have considered ten different models of BRWs, which are presented in Table~\ref{table1}. For simplicity, in all models, the intensity of sources are the identical constants or i.i.d. random variables. Models 1 and 3 represent supercritical BRWs in non-homogeneous and homogeneous non-random media on $\mathbb{Z}$. In these models we put the dying intensity is equal to $1$ and the splitting intensity is equal to $2$. In such models, the first moment experiences exponential growth. Models 1 and 3 are paired with models 2 and 4, respectively. In models 2 and 4 the sources are located in the same way as in models 1 and 3, but intensities are random variables distributed according to the Weibull distribution. Random dying and splitting intensities are selected so that their averages are equal to the intensities of a non-random medium: $\mathbb{E} \left(\text{Weib}(2,1. 13) \right) \approx 1 $ and $\mathbb{E} \left (\text{Weib} (2,2.26) \right) \approx 2 $. In addition, both the random variables under consideration have Gaussian right tail. Their difference has Gaussian right tail too, which was examined numerically. Note that for random intensities, the concept of ``critical'' values does not exist. However, for convenience, from now on we call the BRW in random media ``supercritical'' if for every $x$: $\mathbb{E}\xi_+(x)>\mathbb{E}\xi_-(x) $ and ``critical'' if for every $x$: $\mathbb{E}\xi_+(x) = \mathbb{E}\xi_-(x)$. Model 5 is a critical BRW in a non-random homogeneous medium paired with model 6, a ``critical'' BRW in a random medium. Using models 5 and 6, we will show the difference in the concept of ``critical process'' depending on the randomness of the medium. Finally, models 7--10 are supercritical random and non-random BRWs on $\mathbb{Z}^\mathrm{3}$. The first task of these models is to demonstrate that the effects we are interested in do not depend on the dimension of the integer lattice. The second task is to demonstrate the relative influence of the source configuration on the growth rate of the moments of a total number of particles. The theoretical results for the case of several sources are given, e.g. in~\cite{Y:13,Y:18}.

\begin{table}[h]\small
    \centering
    \begin{tabular}{c c c c c}
 Model number and& Lattice  & Coordinates &  Splitting   & Dying   \\
    brief description & dimension $\mathrm{d}$ &  of sources& intensity $\tilde \xi_+$ & intensity $\tilde \xi_-$\\
 \hline
  \makecell[l]{1. Non-random, supercritical,\\ \hspace{10pt} non-homogeneous} & $1$ & $x = 0$ & 2 & 1 \\
 \makecell[l]{2. Random, ``supercritical'', \\ \hspace{10pt} non-homogeneous} & $1$ & $x = 0$  & $\text{Weib}(2,2.26)$ & $\text{Weib}(2,1.13)$ \\

 \hline
 \makecell[l]{3. Non-random, supercritical, \\ \hspace{10pt}  homogeneous} & $1$ & every point  & 2 & 1 \\
 \makecell[l]{4. Random, ``supercritical'', \\ \hspace{10pt}  homogeneous}& $1$ & every point  & $\text{Weib}(2,2.26)$ & $\text{Weib}(2,1.13)$ \\
 \hline

 \makecell[l]{5. Non-random, critical, \\  \hspace{10pt} homogeneous} & $1$ & every point  & 1 & 1 \\
 \makecell[l]{6. Random,critical,\\ \hspace{10pt}  homogeneous} & $1$ & every point  & $\text{Weib}(2,1.13)$ & $\text{Weib}(2,1.13)$ \\
 \hline

  \makecell[l]{7. Non-random, supercritical, \\ \hspace{10pt}
  a simplex with the side of $\sqrt2$}& $3$ &\makecell{ (1,0,0), (0,1,0),\\(0,0,1)}      & 2 & 1 \\

  \makecell[l]{8. Non-random, supercritical,\\ \hspace{10pt}
   a simplex with the side of $2\sqrt2$} & $3$ & \makecell{ (2,0,0), (0,2,0),\\(0,0,2)}    & 2 & 1 \\

  \makecell[l]{9. Random, ``supercritical'',\\ \hspace{10pt}
   a simplex with the side of $\sqrt2$} & $3$ & \makecell{ (1,0,0), (0,1,0),\\(0,0,1)}   & $\text{Weib}(2,2.26)$ & $\text{Weib}(2,1.13)$ \\

  \makecell[l]{10. Random, ``supercritical'',\\ \hspace{10pt}
   a simplex with the side of  $2\sqrt2$} & $3$ & \makecell{ (2,0,0), (0,2,0),\\(0,0,2)}   & $\text{Weib}(2,2.26)$ & $\text{Weib}(2,1.13)$ \\

\end{tabular}
    \caption{Models under consideration.}\label{table1}
\end{table}

\subsection{Results and discussion}

\begin{figure}[h]
    \centering
     \includegraphics[width=.9\textwidth]{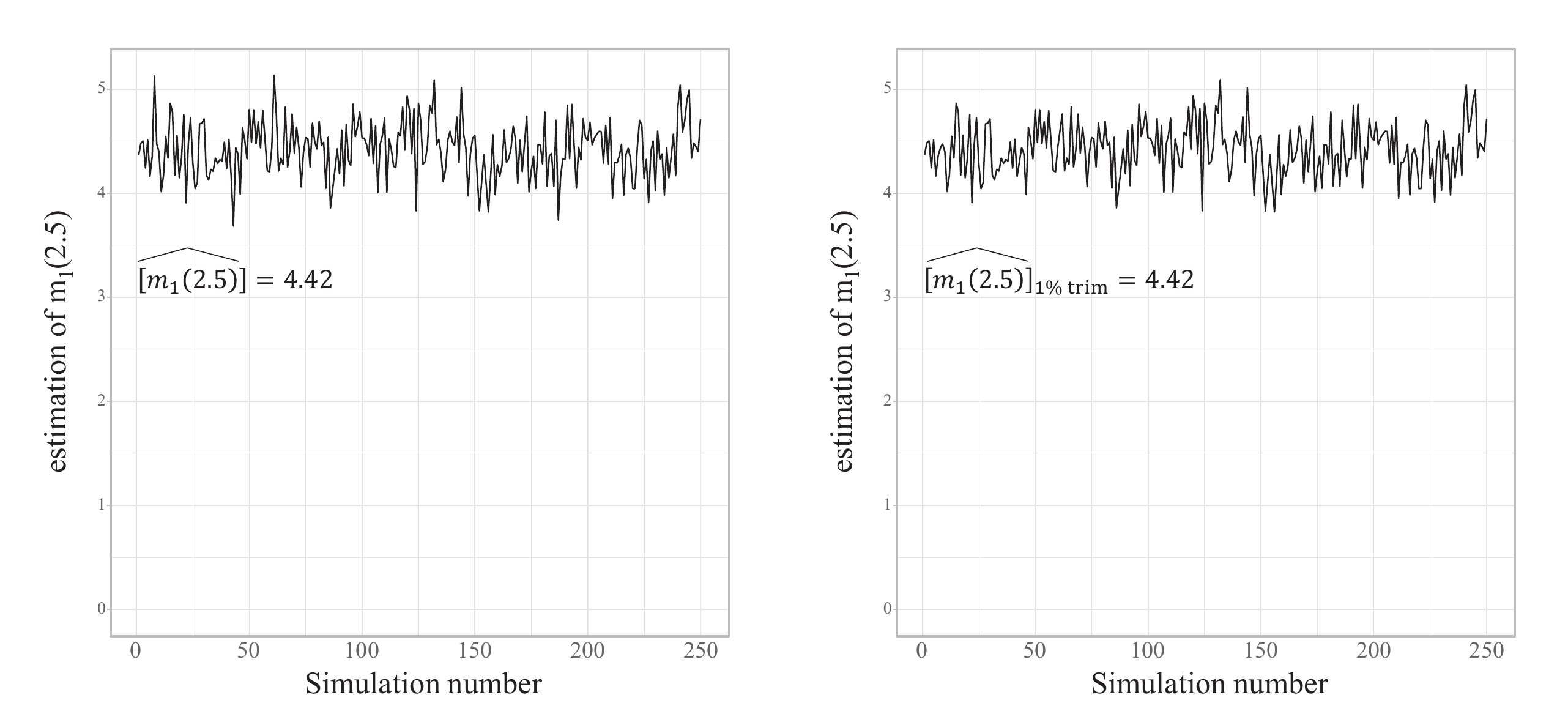}%
    \caption{Estimates and trimmed estimates (right panel) of the first moment of the total number of particles in non-homogeneous non-random BRW.}
    \label{fig1}
\end{figure}

\begin{figure}[h]
    \centering
     \includegraphics[width=.9\textwidth]{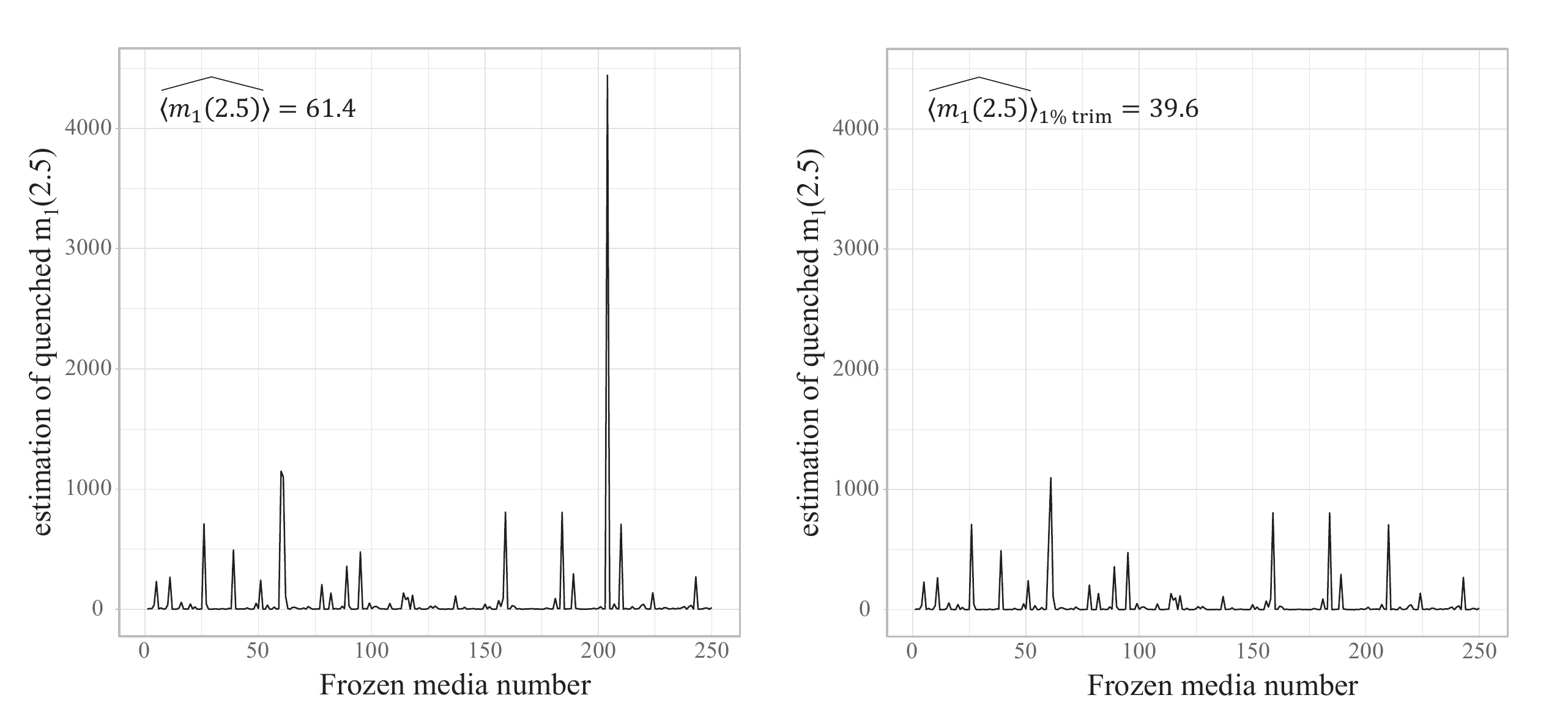}%
    \caption{Estimates and trimmed estimates (right panel) of the first moment of the total number of particles in non-homogeneous random BRW.}
    \label{fig2}
\end{figure}

\begin{figure}[h]
    \centering
     \includegraphics[width=.9\textwidth]{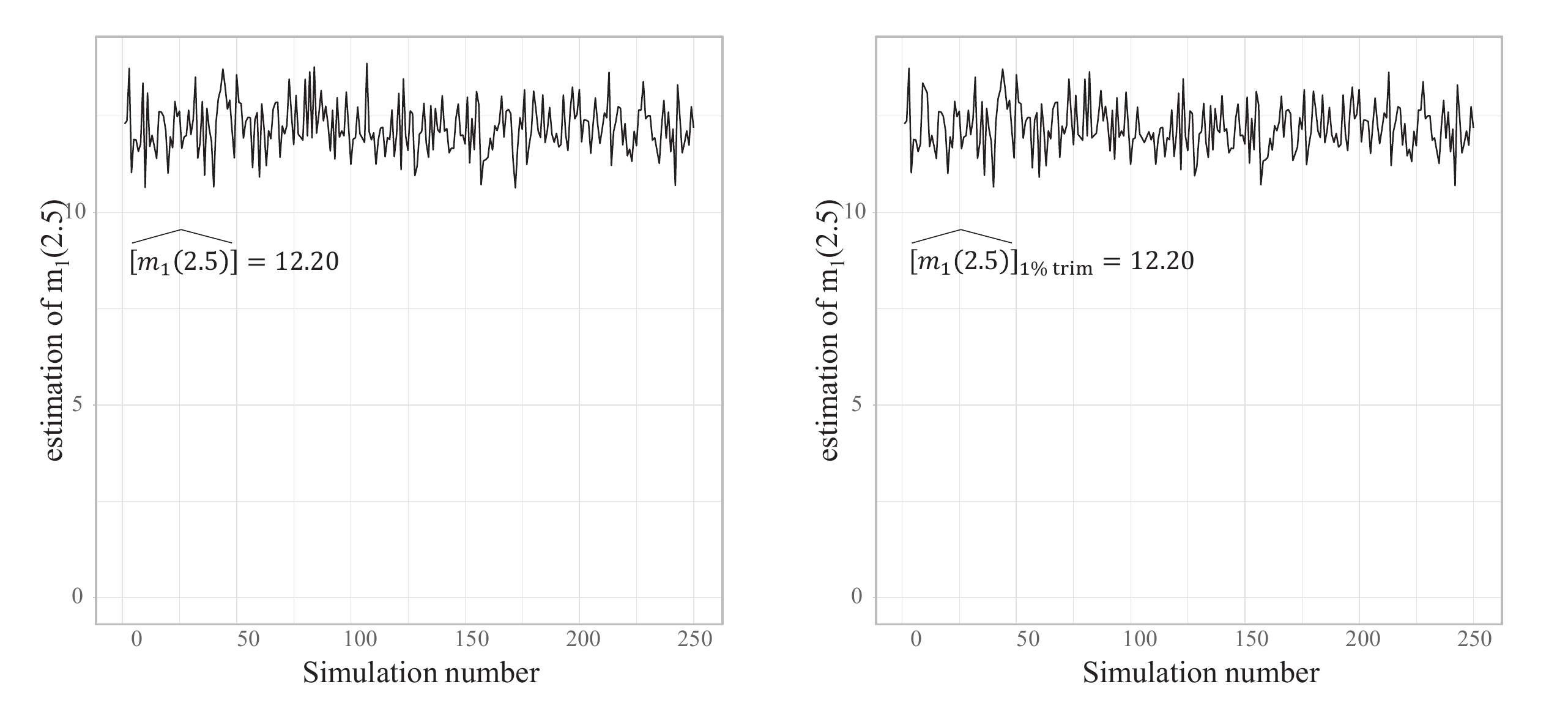}%
    \caption{Estimates and trimmed estimates (right panel) of the first moment of the total number of particles in homogeneous non-random BRW.}
    \label{fig3}
\end{figure}

\begin{figure}[h]
    \centering
     \includegraphics[width=.9\textwidth]{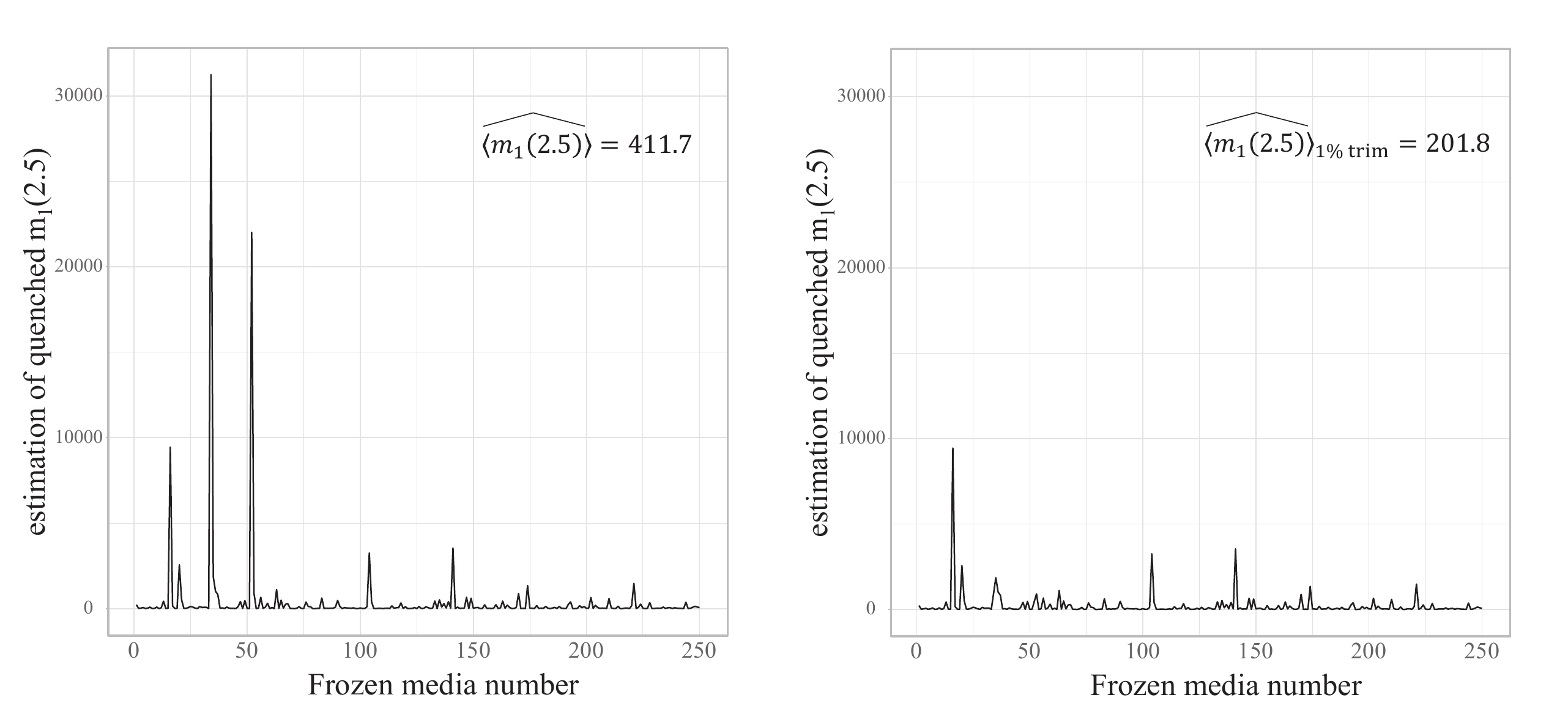}%
    \caption{Estimates and trimmed estimates (right panel) of the first moment of the total number of particles in homogeneous random BRW.}
    \label{fig4}
\end{figure}

Let us start with a graphical representation of the results. In Figures \ref{fig1} and \ref{fig2}, the estimates of the first moment of the total number of particles for model 1 and 2 at time $t = 2.5$ are shown. The moment $t = 2.5$ is chosen due to the visually noticeable demonstration of intermittency. At $t>2.5$, the graphical representation of the quenched moments in random medium becomes uninformative due to high peaks. The behavior of the estimates of quenched moments differs sharply in random and non-random cases. In a non-random medium, the estimates of the first moment look similar to the realizations of a normally distributed random variable. We confirmed the normality by the Shapiro-Wilk test ($p = 0.17$). As we mentioned in section 4.1, this result was expected and confirmed the correctness of the implementation of the Monte Carlo method. In the case of a random medium, the distribution of quenched moments contains rare large-size realizations. We assume that this is a manifestation of the intermittency of the field of quenched moments.  In Figures \ref{fig3} and \ref{fig4}, the similar estimates for models in homogeneous media are presented. These figures confirm the presence of intermittency at finite times.

However, the graphical analysis is not very convenient, so it was necessary to obtain some numerical estimate of the intermittency, particularly for $t = 10$. To do this, we used the ``large fluctuations'' interpretation of the intermittency described in section 3. It states that intermittency leads to the impossibility of a good description of a random variable by its moments. Namely, the main contribution to the moments of a random variable are made by a small number of rare events. Our case implies that the main contribution to the first annealed moment is made by a low number of large first quenched moments.
We decided to estimate this influence by considering the ``$1 \%$ trimmed moments.'' The $1 \%$ trimmed moment value is defined as the average assessed on a sample without lowest and highest $1 \%$ of observations. One percent trimming is equivalent to removing the two largest and two smallest of the $250$ quenched moments and then averaging.  The trimmed samples and the corresponding trimmed moments are shown in the figures \ref{fig1}--\ref{fig4} on the right panels. The trimmed average in the case of a non-homogeneous non-random medium decreased by $36 \%$ compared to its original value. In contrast, the same value in a random medium changed only in the second decimal place or $0.00 \%$. A similar result is observed for homogeneous media. Therefore we need to emphasize that intermittency occurs not only in a homogeneous medium, but also in a non-homogeneous medium with only one source of branching.

Let us  introduce the function $R(t)$ as a measure of intermittency which estimates how much the first annealed moment depends at time $t$ on rare realizations of the field:

\[R(t) =
\begin{cases}

\left. \widehat{\left[m_1(t) \right]} \middle/  \widehat{\left[m_1(t) \right]}_{\text{$1 \%$ trim}}
\right. \text{in the case of non-random media}, \\
\left. \widehat{\langle m_1(t) \rangle} \middle/ \widehat{\langle m_1(t) \rangle}_{\text{$1 \%$ trim}} \right. \text{in the case of random media}.

\end{cases}
\]

\begin{figure}[h]
    \centering
     \includegraphics[width=.9\textwidth]{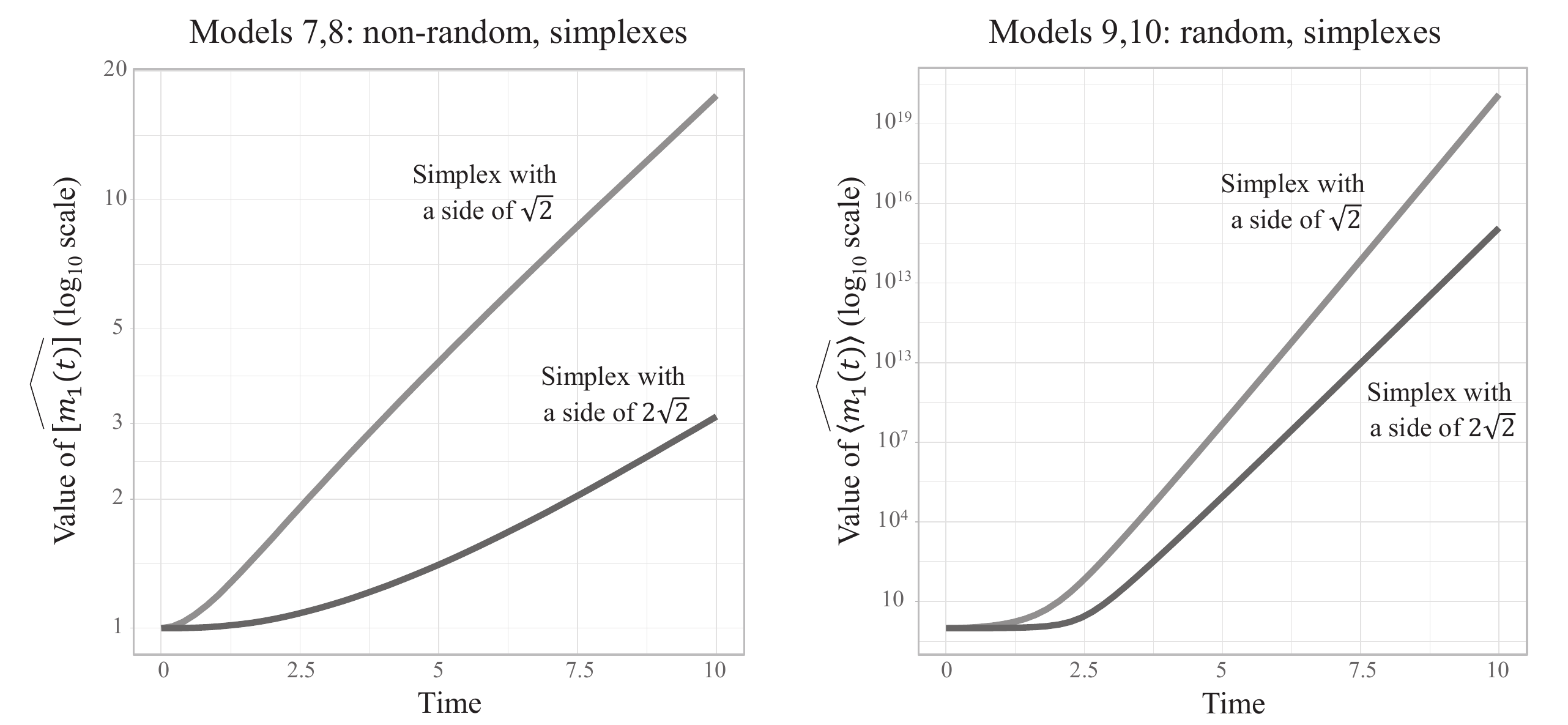}%
    \caption{Moments and annealed moments in models 7--10}
    \label{fig5}
\end{figure}

We estimated $R(10)$ for all models from Table~\ref{table1}, and the results are shown in Table~\ref{table2}. As one can see, for models with a non-random medium, $R(10)$ takes values equal to $1$ with the third decimal place accuracy. In turn, this ratio is significantly greater than $1$ for models with random medium. Thus, for all dimensions and source configurations, the assumption of randomness radically changes the behavior of the BRW.   When comparing the pair of models 1--2 with the pair 3--4, one can observe that a homogeneous medium in the supercritical case generates a greater increase in moments. This result can be explained as follows. In a homogeneous medium, the random walk only slightly delays branching by moving the particle from one branching center to another. While in a non-homogeneous medium, the walk can significantly delay branching, forcing the particle to wander outside the source. This effect is especially important in spaces of dimension greater than 3, where a simple symmetric random walk can ``remove'' a point from the source of branching permanently. This effect is discussed in more detail in~\cite{YarBRW:r}. On the example of models 5 and 6, a theoretical result on the concept of ``critical values'' of the BRW is illustrated. These values is the essential characteristic of the process in a non-random medium, and, at the same time, they does not appear in any way in a non-random medium. Finally, models 7--10 show that the growth rate of moments also depends on the configuration of sources. For clarity, we have presented estimates of moments for the models 7--10 in Figure~\ref{fig5}. It can be seen that the configuration determines the exponent power in the asymptotic of the first moment in both random and non-random media. Moreover, the influence of the configuration seems obvious: the farther the sources are from each other, the slower the growth.

\begin{table}[h]\small
    \centering
    \begin{tabular}{c c c c}
 Model number and&   $\widehat{\left[m_1(10) \right]} $ or $\widehat{\langle m_1(10) \rangle}$&
                    $\widehat{\left[m_1(10) \right]}_{\text{$1 \%$ trim}} $ or    & $R(10)$  \\
    brief description &  & $\widehat{\langle m_1(10) \rangle}_{\text{$1 \%$ trim}}$ &  \\
 \hline
 \makecell[l]{1. Non-random, supercritical,\\ \hspace{10pt} non-homogeneous}  & $111.0$ & $110.9$ & $1.00$ \\
 \makecell[l]{2. Random, ``supercritical'', \\ \hspace{10pt} non-homogeneous} & $6.5 \cdot 10^{12}$ & $4.9 \cdot 10^{10}$ & $131$ \\

 \hline
 \makecell[l]{3. Non-random, supercritical, \\ \hspace{10pt} homogeneous}  & $2.5 \cdot 10^4$ & $2.5 \cdot 10^4$ & $1.00$ \\
 \makecell[l]{4. Random, ``supercritical'', \\ \hspace{10pt} homogeneous} & $6.3 \cdot 10^{16}$ & $2.2 \cdot 10^{15}$ & $28$ \\
 \hline

 \makecell[l]{5. Non-random, critical, \\ \hspace{10pt} homogeneous}  & $0.997 $ & $ 0.997 $ & $1.00$ \\
 \makecell[l]{6. Random, ``critical'',\\ \hspace{10pt} homogeneous}  & $ 59.8 $ & $44.7$ & $1.3$ \\
 \hline

 \makecell[l]{7. Non-random, supercritical, \\
  \hspace{10pt} simplex with the side of $\sqrt2$} & $17.5 $ & $17.5$ & $1.00$ \\

 \makecell[l]{8. Non-random, supercritical\\
   \hspace{10pt} simplex with the side of $2\sqrt2$}  & $3.1$ & $3.1$ & $1.00$ \\

 \makecell[l]{9. Random, ``supercritical'',\\
   \hspace{10pt} simplex with the side of $\sqrt2$}  & $1.3 \cdot 10^{20}$ & $2.3 \cdot 10^{13}$ & $5.3 \cdot 10^6$ \\

 \makecell[l]{10. Random, ``supercritical'',\\
   \hspace{10pt} simplex with the side of $2\sqrt2$}  & $1.2 \cdot 10^{15}$ & $2.3 \cdot 10^4$ & $5.1$ \\

\end{tabular}
\caption{Simulation results}\label{table2}
\end{table}

\begin{figure}[h]
    \centering
     \includegraphics[width=.9\textwidth]{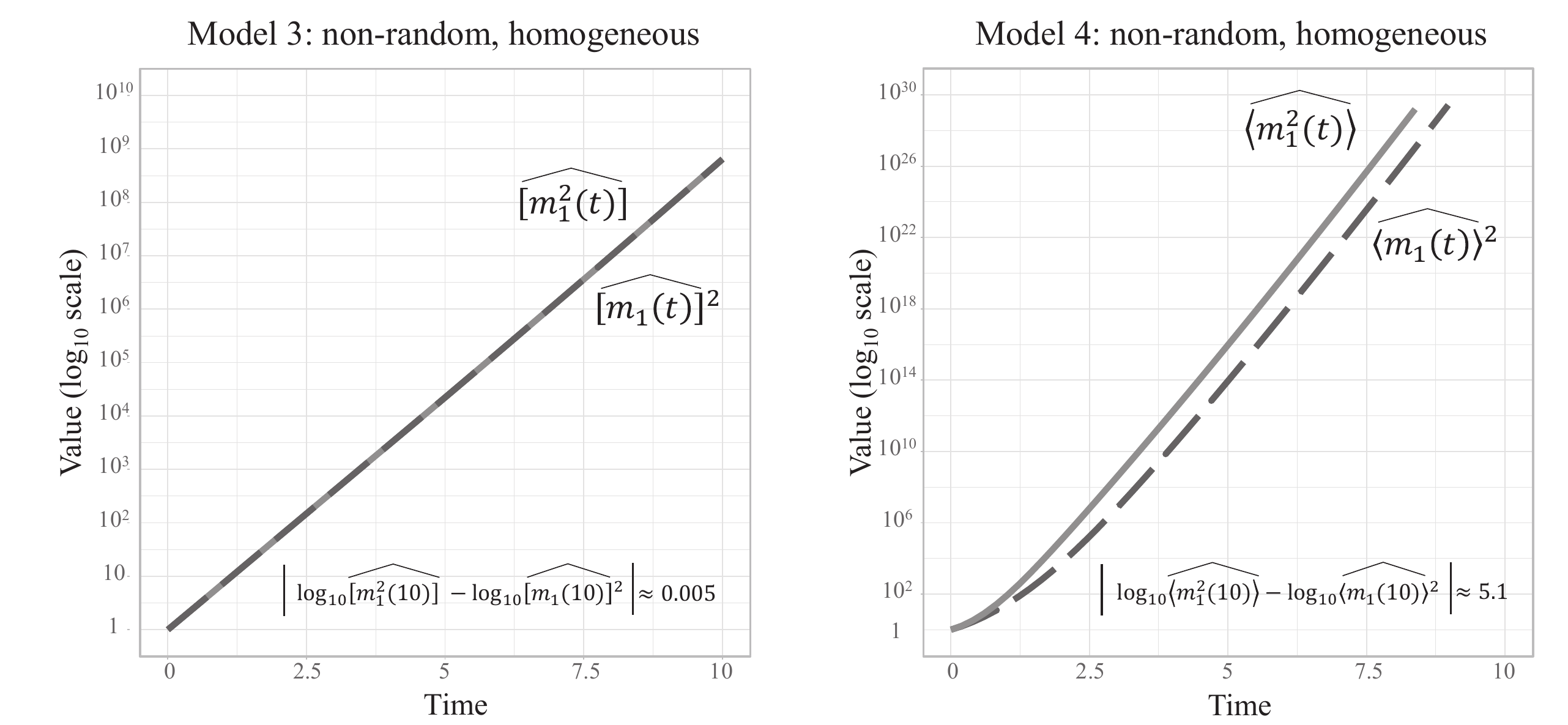}%
    \caption{Comparison of the growth of the first two annealed moments}
    \label{fig6}
\end{figure}

Note that the property of ``high peaks'' of the field of moments is not the only characteristic of intermittency and is a consequence derived from the definition, see Section~\ref{S3}. We are rather interested in the irregularity of the growth of moments, which is equivalent to the definition of intermittency. In particular, for BRW in a random medium, the following property is valid:
$\langle m_{1}(t) \rangle^2 \ll \langle m_1^2 (t) \rangle$, as $t \to \infty $.
We reformulate this property in terms of logarithms: $\log_{10} \langle m_1^2 (t) \rangle - \log_{10} \langle m_{1} (t) \rangle^2 > C$, as $t \to \infty, C \in \mathbb{R}_+ $. Note that in practice, it cannot be shown that $C$ is exactly zero. One can only test the hypothesis that $C = 0$, but it requires generating a sample of estimates of annealed moments, which is computationally difficult. In this paper, we only estimated the irregularity of growth for the obtained estimate of the annealed moment. Figure \ref{fig6} shows graphs for annealed moments in the case of random and non-random media. For a non-random medium, both log-values behave similarly with a difference of less than $0.005$. While for a non-random medium, the constant $C$ stabilizes at approximately $5.1$. The figure and the estimates confirm that we could observe the phenomenon of intermittency qualitatively and quantify the degree of irregularity of the moment growth.

\section{Conclusion}

The intermittency can be observed in random media even over finite time intervals. The Monte Carlo method can be used as an exploratory tool for BRW with arbitrary distributions of branching intensities and field structures. However, in random media, estimates of moments become unstable due to the presence of intermittency.

\textbf{Acknowledgements.} The research was supported by the Russian Foundation for the Basic
Research (RFBR), project No. 20-01-00487.

\end{document}